\mag\magstep1
\def\SetAuthorHead#1{}
\def\SetTitleHead#1{}
\def\NoindentAfter{\everypar={\setbox0=\lastbox\everypar={}}}
\def\H#1\par#2\par{{\baselineskip=15pt\parindent=0pt\parskip=0pt
 \leftskip= 0pt plus.2\hsize\rightskip=0pt plus.2\hsize
 \bf#1\unskip\break\vskip 4pt\rm#2\unskip\break\hrule
 \vskip40pt plus4pt minus4pt}\NoindentAfter}
\def\HH#1\par{{\bigbreak\noindent\bf#1\medbreak}\NoindentAfter}
\def\HHH#1\par{{\bigbreak\noindent\bf#1\unskip.\kern.4em}}
\def\th#1\par{\medbreak\noindent{\bf#1\unskip.\kern.4em}\it}
\def\endth{\medbreak\rm}
\def\pf#1\par{\medbreak\noindent{\it#1\unskip.\kern.4em}}
\def\df#1\par{\medbreak\noindent{\it#1\unskip.\kern.4em}}

\let\rk\df\let\endrk\enddf
\let\Roster\bgroup\let\endRoster\egroup
\def\\{}\def\text#1{\hbox{\rm #1}}
\def\mop#1{\mathop{\rm\vphantom{x}#1}\nolimits}
\def\MaxReferenceTag#1{}
\def\qedbox{\vrule width2mm height2mm\hglue1mm\relax}
\def\qed{\ifmmode\qedbox\else\hglue5mm\unskip\hfill\qedbox\medbreak\fi\rm}

\def\cite#1{{\bf[#1]}}
\def\Em#1{{\it #1\/}}\let\em\Em
\def\Bib#1\par{\bigbreak\bgroup\centerline{#1}\medbreak\parindent30pt
 \parskip2pt\frenchspacing\par}
\def\endBib{\par\egroup}
\newdimen\Overhang
\def\rf#1{\par\noindent\hangafter1\hangindent=\parindent
     \setbox0=\hbox{[#1]}\Overhang\wd0\advance\Overhang.4em\relax
     \ifdim\Overhang>\hangindent\else\Overhang\hangindent\fi
     \hbox to \Overhang{\box0\hss}\ignorespaces}

\def\bbZ{{\Bbb Z}}
\def\Coordinates{\bigbreak\bgroup\parindent=0pt\obeylines}
\def\endCoordinates{\egroup}


\def\Cal#1{{\cal #1}}
\def\hexnumber#1{\ifcase#1 0\or 1\or 2\or 3\or 4\or 5\or 6\or 7\or 8\or
 9\or A\or B\or C\or D\or E\or F\fi}
 
\font\tenmsa=msam10
\font\sevenmsa=msam7
\font\fivemsa=msam5
\newfam\msafam
\textfont\msafam=\tenmsa
\scriptfont\msafam=\sevenmsa
\scriptscriptfont\msafam=\fivemsa

\font\tenmsb=msbm10
\font\sevenmsb=msbm7
\font\fivemsb=msbm5
\newfam\msbfam
\textfont\msbfam=\tenmsb
\scriptfont\msbfam=\sevenmsb
\scriptscriptfont\msbfam=\fivemsb

\def\Bbb{\fam\msbfam\tenmsb}
\vsize 20cm
\hsize 14cm
\def\rtimes{\times \hskip -0.6em \vert}

 
\def \ppp{{\cal P}}
\def \atp{AT$\ppp$}
\def \Aut{\mop {Aut}}
\let \eval=\Overline
\def\mto{{\buildrel M \over \to}}
\def\mstar{{\buildrel M *\over \to}}
\def\NR{NEED REFERENCE}

\def\today{\number\day\space\ifcase\month\or January\or 
February\or March\or April\or May\or June\or July\or 
August\or September\or October\or November\or December
\fi\space\number\year}
\font\small=cmr8
\font\smallit=cmti8
\font\smallbf=cmbx8
\long\def\omit#1{\relax}
\omit{\footline={\par \vbox{\leftskip = 0pt\rightskip = 0pt  
\hrule\vglue-0.04truein
\vskip 0.2cm\leftskip = 0pt {\small PPP groups}
 \hskip 2cm  {\smallbf page  \folio }
\hfill  {\smallit printed \today} }\par }}

\input amssym.def 
\input amssym.tex

\SetTitleHead{}
\SetAuthorHead{}

\H Parallel Poly Pushdown Groups 

Gilbert Baumslag, Michael Shapiro and Hamish Short

\HH 1. Introduction 

The theory of automatic groups has attracted wide attention for the
last ten years or so and brought together a large and varied collection
of groups. These groups admit a simple description in terms of
finite state automata - we refer the reader to {\S 2} for details
and the various terms that we shall make use of in this
introduction.

One of the original motivations for this theory was to provide a
means for computing in the fundamental groups of geometric
3-manifolds. The goal was, in part, accomplished by \cite{ECHLPT}
and \cite{S} when they, independently, prove the following theorem:
if M is a 3-manifold that obeys Thurston's geometrization conjecture,
then $\pi_1(M)$ is automatic if and only if M does not contain a
closed Nil or Sol manifold in its connected sum decomposition.

One generalization of the class of automatic groups is the class of
asynchronously automatic groups. While this class successfully
captures such non-automatic groups such as the Baumslag-Solitar
groups \cite{BGSS}, \cite{ECHLPT}, it still does not include the fundamental
groups of Nil \cite{ECHLPT} and Sol \cite{Br} manifolds.

Bridson has shown that the fundamental groups of Nil and Sol 
manifolds have asynchronous combings.  That is to say, one can find a 
normal form in such a group with the so-called asynchronous fellow 
traveler property.  Since the groups in question are not 
asynchronously automatic, the language of such a normal form is not 
the language of a finite state automaton.  Bridson and Gilman 
\cite{BGi} have investigated the computational complexity of these 
asynchronous combings.  They show that they can not be discovered by 
means of a pushdown automaton, but they can be discovered by means of 
a nested stack automaton.

Here we shall take a slightly different approach. Motivated by the
notion of computations being carried out in parallel, we allow for
the use of as many pushdown automata as we need, in order to describe
a variety of different groups. This gives rise to to a class of
languages that we call {\em{parallel poly-pushdown}}, allowing us
to generalize automatic groups to {\em {parallel poly-pushdown groups}}.
Passing from finite state automata to pushdown automata increases
the space costs from a bounded amount of memory to a linearly
bounded amount of memory. However this does not increase time costs.
Each machine still processes its input in linear time. Thus, from
a computational point of view, little is lost.

The resulting class of parallel poly-pushdown groups,
which we shall denote here simply by $\Cal P$, has a number of
properties in common with automatic groups. For example, we shall
prove that both the free and the direct product of two groups
in $\Cal P$ is again in $\Cal P$. In addition, we shall prove
that such parallel poly-pushdown groups are recursively presentable.
In fact even more is true:

\th Theorem 4.2

$\Cal P$-groups have solvable word problem. \endth

The similarity between the two classes ends here. One of the
most interesting points of departure, which underlines the
significance of Theorem 4.2, is a consequence of the following

\th Theorem 5.4

The standard wreath product $U \wr T$ of the parallel poly-pushdown 
group U by the parallel poly-pushdown group T is again a parallel 
poly-pushdown group.
\endth

Now the wreath product $W$ of a finitely presented group $U$
by a finitely presented group $T$ is finitely presented if and
only if either $U=1$ or $T$ is finite \cite{Ba}. Since the infinite
cyclic group $C$ belongs to $\Cal P$, it follows that the
finitely generated but not finitely presented group $C \wr C$
belongs to $\Cal P$, i.e., $\Cal P$ is not contained in the class of
finitely presented groups, in direct contrast to the class of
automatic groups.

Now a finite generated nilpotent group is automatic if and only
it is virtually abelian \cite{ECHLPT}. Again we have a pronounced
difference between the two classes of groups:

\th Theorem 5.2

Let A be a finitely generated central subgroup of the group G.
If G/A is automatic, then $G \in \Cal P$.
\endth

It follows thence that

\th Corollary 5.3

Every finitely generated nilpotent group of class at most
two is parallel poly-pushdown.
\endth

We have been unable to decide whether every finitely generated
nilpotent group is parallel poly-pushdown. However
we have been able to prove

\th Theorem 5.6

Suppose that $H$ is a $\ppp$ group, and $\varphi:H \to \Aut {\Bbb
Z}^n$.  Then ${\Bbb Z}^n \rtimes_\varphi H$ is $\ppp$. \endth

It follows that there are parallel poly-pushdown groups
which are nilpotent of arbitrary class. Indeed, every finitely 
generated torsion-free nilpotent group embeds in a parallel-poly 
pushdown group, as we see from

\th Theorem 6.1

Let $n$ be any positive integer. Then the group of all integral upper
uni-triangular matrices of degree $n$, is parallel poly-pushdown.
\endth

This implies that every finitely generated torsion-free nilpotent
group is a subgroup of a parallel poly-pushdown group,
namely the above (nilpotent) group of matrices.

Finally we remark (Corollary 5.10) that $\Cal P$ 
contains the fundamental groups of
all 3-manifolds which obey Thurston's geometrization conjecture.

\HH 2. Preliminaries

Finite state automata and regular languages are now fairly well known 
to geometric group theorists.  
We record the basic definitions for completeness.

Given a finite set $A=\{a_{1},\dots,a_{k}\}$ the free monoid on $A$ is
denoted $A^{*}$.
Thus $A^{*}$ consists of  strings 
$w=a_{i_{1}}\dots a_{i_{n}}$ where $a_{i_{j}}\in A$ and $n\ge 0$.  
The elements of $A$ are called the \Em{letters} of the  
\Em{alphabet} $A$, and the 
elements of $A^{*}$ are called \Em{words}.
The \Em{length} of the word $w=a_{i_{1}}\dots a_{i_{n}}$ is $n$,
written $\ell(w)=n$.  
If $n=0$ then $w$ is the \Em{empty word} denoted by $\epsilon$.  
A subset of $A^{*}$ is called a \Em{language}.

A finite state automaton is an idealization of a machine which has an
input tape and finite amount of internal memory (and hence finitely
many states that memory can be in).  
The automaton reads the input tape and changes its  state 
according to its current state, together with
the letter just read from the tape. 
More formally, a \Em{finite state automaton} over the alphabet 
$A$ is a quintuple
$(S,A,\tau,s_{0},Y)$ where $S$ is a finite set of \Em{states},
$s_{0}\in S$ is the \Em{start state}, $\tau:S\times A \to S$ is the
\Em{transition function} and $Y\subset S$ is the set of \Em{accept
states}.  
A finite state automaton $M$ determines a language
$L(M)\subset A^{*}$ in the following way.  
For each word $w=a_{i_{1}}\dots a_{i_{n}} \in A^{*}$, 
let $t_{0}=s_{0}$ and for $1\le j\le n$ let
take $t_{j}=\tau(t_{j-1}, a_{i_{j}})$.
Now let 
$$L(M) = \{
w=a_{i_{1}}\dots a_{i_{n}}\in A^{*} \mid n\ge 0, t_n \in Y\}.$$ 
A language is \Em{regular} if it is the language of some finite state
automaton.

The \Em{concatenation} of two words $u$ and $v$ in $A^{*}$ is their 
product $uv$.  
Given two languages $L\subset A^{*}$ and $M\subset A^{*}$, 
the concatenation of these languages is
$$LM=\{uv \mid u\in L, \ v \in M\}.$$
The Kleene star of a language $L$ defined to be 
$$L^{*}= \{\epsilon\}\cup L \cup LL \cup LLL \dots $$
(recall that $\epsilon$ denotes the empty word.)
The class of regular languages is closed under the set operations of union, 
intersection and complementation and under the language theoretic 
operations of concatenation and Kleene star.   In fact the class of 
regular languages is the smallest class of languages which contains 
all finite languages and is closed under these operations.

Pushdown automata are somewhat less common in geometric group theory.
Roughly speaking, a pushdown automaton is a finite state automaton
which also has control of a stack.  It can write to the top of the
stack and erase from the top of the stack.  At any point in a
computation it chooses its next move based on its current internal
state, the top letter of the stack and possibly a letter read from its
input tape.  
More formally, a \Em{deterministic pushdown automaton} is
a 7-tuple $(S,A,B,\tau,s_{0},Z_{0},Y)$ where
\item{1.} $S$ is a finite set of \Em{states}.
\item{2.} $A$ is a finite \Em{input alphabet}.
\item{3.} $B$ is a finite \Em{stack alphabet}.
\item{4.} The transition function $\tau$ is a map from a 
subset $\Sigma$ of $(A\cup\{\epsilon\})\times B\times S$ to $S\times 
B^{*}$.
\item{5.} $s_{0} \in S$ is the \Em{start state}.
\item{6.} $Z_{0} \in B$ is the \Em{start symbol}.
\item{7.} $Y\subset S$ is the set of \Em{final states}.

Further conditions on the set $\Sigma$  will be given below.

We now describe how the automaton 
$M= (S,A,B,\tau,s_{0},Z_{0},Y)$ accepts or rejects a word $w\in A^{*}$.  
At any point in a computation, the current unread content of the
input tape is a word $u\in A^{*}$, 
the current state is an element $s\in S$ 
and the current content of the stack is a word $z \in B^{*}$.  
The \Em{instantaneous description} of this state of affairs is the triple
$(u,z,s)$.  
Suppose that $u=av$ with $a \in A \cup \epsilon$ and $z=yb$ with $b\in B$.  
Then $M$ makes the transition  $(av,yb,s) \mto (v, y\alpha,t)$ 
where $(t,\alpha)=\tau(a,b,s)$.  
That is, $M$ reads the first letter (or no letter) of its remaining input, 
replaces the top letter of its stack with a word 
(possibly the empty word, possibly the letter just erased) 
and makes a transition to a new internal state.  
To ensure that the operation of $M$ is deterministic,
i.e., that $M$ has exactly one possible transition for
each instantaneous description, the
following requirement is imposed on $\Sigma$. 
Suppose that $b$ is a stack letter and $s$ is a state.  
Then either 
(i) $(\epsilon,b,s)\in \Sigma$, and $(a,b,s)\notin \Sigma$ for $a\in A$, 
or (ii) $(\epsilon,b,s)\notin \Sigma$, and $(a,b,s)\in \Sigma$ for
all $a\in A$.  
 
There are several ways to define what it means for $M$ to accept a
word. Let $\mstar$ be the reflexive and transitive closure of $\mto$.  
We say $M$ \Em{accepts $w$ by empty stack} if
$(w,Z_0,s_0)\mstar(\epsilon,\epsilon,s)$.  
We say $M$ \Em{accepts $w$ by final state} 
if $(w,Z_0,s_0)\mstar(\epsilon,z,s)$ with $s \in Y$.  These two
notions are equivalent in the following sense.

\th Proposition 2.1

Let $L\subset A^*$.  Then there is a deterministic pushdown automaton
$M$ such that $L$ is the set of words which $M$ accepts by empty
stack if and only there is a deterministic pushdown automaton $M'$ so
that $L$ is the set of words which $M'$ accepts by final state. 
\endth

If $L$ is a language determined by a deterministic pushdown automaton 
in either of these ways, 
then $L$ is said to be \Em{the language of a deterministic pushdown automaton}
(or $L$ is a deterministic context--free language).  
In view of Proposition 2.1, the automata constructed in this article 
will use either acceptance by empty stack or acceptance by final state, 
as seems most convenient. 
We shall usually not provide an explicit
construction of the states, transition function, etc.  
However, we hope that our descriptions will be sufficiently clear 
in each case for the dedicated reader to be able to construct explicitly 
the corresponding automaton as a straightforward (if tedious) exercise.

The membership problem for deterministic pushdown languages can be solved in an 
efficient manner in the following sense:

\th Proposition 2.2

Let $L\subset A^*$ be the language accepted by a deterministic
pushdown automaton.

There is an algorithm to decide whether or not words $w\in A^{*}$
lie in $L$ taking time proportional to the length of the 
word $w$ considered.
\endth

When $L$ is the language of a general non--deterministic pushdown automaton,
then one can prove that there exists such an algorithm taking time 
proportional to the cube of the length of $w$.
For proofs of these results, and a general treatment of pushdown automata
(context--free grammars etc.), see for instance \cite{HU} 
(especially chapters 5,6 and 10).

One simple use of a pushdown stack is as a counter.  We can use the 
stack to keep track of an element of $\bbZ$ in the following manner.  
The number $n\ge 0$ might be represented by placing $n$ copies of some 
symbol (say $+1$) on the stack.  The number $-n \le 0$ might be 
represented by placing $n$ copies of some symbol (say $-1$) on the 
stack.  If we choose some bound $C$ in advance then we can build a 
pushdown automaton $M$ which is capable of adding any number $m$, $-C 
\le m \le C$ to the current contents of the stack.  If $m$ and the 
stack contents $n$ have the same sign, $M$ pushes $|m|$ copies of the 
appropriate symbol onto the stack.  If $m$ and the stack contents $n$ 
differ in sign, $M$ begins by popping from the stack.  It continues 
until it has either done this $|m|$ times or has exhausted the stack.  
If it exhausts the stack before popping $|m|$ times, it continues by 
pushing the appropriate number of symbols of the opposite sign onto 
the stack.  Since $|m|$ is bounded, all of this can be programmed into 
the states of $M$ to be performed upon consuming a single letter of 
input.

Similarly, $M$ can be designed so that if the number $n$ is contained 
in the stack and $M$ encounters a string $x^{m}$ as input, then $M$ 
can consume the string $x^{m}$ and add $m$ to the stack.

The languages studied in this article are the intersections of
finitely many deterministic pushdown languages, so the same efficiency 
result holds.  

We make the following definition:

\df Definition 

A language $L \subset A^*$ is 
said to be \Em{parallel poly pushdown (${\cal P}$)} 
if there are finitely many languages  $L_i$, $i=1,\ldots,k$ of deterministic 
pushdown automata, such that $L= \bigcap_{i=1}^k L_i$.
 
 Notice that the class of deterministic  pushdown languages is not 
closed under the operation of intersection, as is shown by the
example considered in \cite{HU}, \S 6.2:
$$\{a^ib^ic^i\mid i>0\} = 
\{a^ib^ic^j\mid i,j>0\}\cap\{a^jb^ic^i\mid i,j>0\} \ .$$

Our aim is to use the class of $\ppp$ languages to describe 
normal forms in certain groups.
To describe multiplication in the groups concerned, 
following the ideas used in the theory of automatic groups,
we describe a method for deciding whether or not two given 
words in normal form correspond to group elements which differ by 
multiplication on the right by a generator. 
This leads to the notion  of
an asynchronous two tape parallel poly pushdown language.  

Given an alphabet $A$, let $A^\sharp$
and $A^\flat$ be disjoint sets of the same cardinality as $A$,
with isomorphisms 
$u:A^*\to A^{\sharp*}$ and $v:A^*\to A^{\flat*}$ 
denoted by $u \mapsto u^\sharp$ and $v \mapsto v^\flat$.
Given a pair of words $(u,v) \in A^* \times A^*$, a \Em{shuffle},
$\sigma(u,v)$ is defined to be a word $u_1^\sharp v_1^\flat \ldots u_k^\sharp
v_k^\flat \in (A^\sharp \coprod A^\flat)^*$ 
so that $u = u_1 \ldots u_k$, $v=v_1 \ldots v_k$ 
and   $u_i\neq\epsilon$ if $i>1$, and $v_j\neq\epsilon$ if $j<k$.  
In spite of the notation $\sigma$ is not a function.  

Let $\$ $ be a symbol which is not in the alphabet $A$.  
An \Em{asynchronous two tape pushdown automaton} over the alphabet $A$ 
is a deterministic pushdown automaton whose input alphabet is 
$(A\cup\{\$\})^\sharp \coprod (A\cup\{\$\})^\flat$,
and whose states are partitioned into two disjoint sets
$S_\sharp$ and $S_\flat$.  
Given an asynchronous two tape pushdown automaton $T$ over the alphabet $A$,
and a pair of words $(u,v) \in A^*\times A^*$, 
we say that $(u,v)$ is \Em{accepted} by $T$
if there is a shuffle $\sigma=\sigma(u\$,v\$)$ so that $\sigma$ is
accepted by $T$ and when reading $\sigma$, $T$ reads only
$(A \cup\{\$\})^\sharp$ letters while in $S_\sharp$ states and only
$(A \cup\{\$\})^\flat$ letters while in $S_\flat$ states.  
We call such a $\sigma$ an \Em{acceptance} of $(u,v)$.  
The deterministic nature of the automaton $T$ has the following 
consequence for the efficiency of $T$:

\th Lemma 2.3

For each pair $(u,v)\in A^{*}\times A^{*}$, there is at most one 
shuffle $\sigma=\sigma(u\$,v\$)$ which is accepted by $T$.  Moreover, 
there is an algorithm taking time proportional to the sum of the 
lengths of the words $u,v$, which determines whether or not such a 
shuffle exists.\qed

\endth
 
We say that the \Em{language} of $T$ is the set of pairs accepted by 
$T$.  Strictly speaking $T$ cannot read an $A^{\sharp}$ letter while 
in a $S_{\flat}$ state, nor can it read an $A^{\flat}$ letter while 
in a $S_{\sharp}$.  Speaking colloquially, we will say that $T$ 
``goes to a fail state'' in such a situation.

We say that a subset of $A^* \times A^*$ is an \Em{asynchronous two
tape parallel poly pushdown} language (AT${\cal P}$) if it is the
intersection of the languages of finitely many asynchronous two tape
pushdown automata.  We emphasize that these machines are taken to be 
deterministic.  
It follows immediately from Lemma 2.3 that membership of such a 
language can be determined by an algorithm taking time proportional
to the sum of the lengths of the pair of words involved.

It seems likely that if we pass to nondeterministic machines the 
analogous algorithms could take exponential time.

\HH 3. Elementary properties of $\ppp$ languages

\th Proposition 3.1

Suppose that $L$ and $M$ are $\ppp$, and that $R$ is regular.  Then
the following are  ${\cal P}$:

{\it 1.}  $L\cap M$.

{\it 2.} $L \cup R$.

{\it 3.} $L - R$.

{\it 4.} If $M$ and $L$ are $\ppp$ languages over disjoint 
alphabets, and $\epsilon\not\in L\cup M$, then $LM$ and $(LM)^{*}$ are $\ppp$.
  
If $M$ and $L$ are $\ppp$ languages over disjoint alphabets, 
and $\epsilon \in L\cap M$, 
then  $$M(L-\{\epsilon\}M-\{\epsilon\})^{*}L$$ is $\ppp$.
\endth

\pf  Proof

{\it 1.} We have $L= \bigcap_{i=1}^m L_i$, $M= \bigcap_{i=1}^n
M_i$, so $L\cap M =  (\bigcap_{i=1}^m L_i) \cap (\bigcap_{i=1}^n
M_i)$.

{\it 2.}  We have $L \cup R = \bigcap_{i=1}^m (L_i \cup R)$.
But the union of a deterministic pushdown language with a regular
language is itself a deterministic pushdown language.  (See \cite{HU}.)

{\it 3.} $L-R = L \cap (A^* - R)$.  But $A^* -R$ is regular,
hence deterministic pushdown, so we are done. 

{\it 4.} We have $L= \bigcap_{i=1}^m L_i$, $M= \bigcap_{i=1}^n M_i$, 
where  $L_i$  and $M_i$ are the languages of deterministic
pushdown automata $P_i$ and $Q_i$ 
(we suppose that these automata accept by final state).  
Notice that $LM = \cap_{i,j}L_iM_j$.  
Thus it suffices to construct a deterministic pushdown automaton $S_{ij}$ 
which accepts the language $L_iM_{j}$.  
$S_{ij}$ starts in the start state of $P_i$ and
continues with the operation of $P_i$ until it encounters a generator
from the alphabet of $M$.
If it is not in an accept state of $P_i$, it goes to a fail state and
reads the rest of the tape.  
If it is in an accept state of $P_i$, it empties the stack, and goes 
into the start state of $Q_j$.  
The automaton $S_{ij}$ accepts if and only if it ends in
an accept state of $Q_j$.

The same construction generalizes to $(LM)^{*}$, as the conditions 
given on the languages ensure that
$\bigl((\bigcap_{i=1}^m L_i)(\bigcap_{j=1}^n M_j)\bigr)^{*} = 
\bigcap(L_{i}M_{j})^{*}$. 

For the final more complicated case,
we suppose that the there is no $\epsilon$ transition from the
initial states of the machines for the languages used to
define $L$ and $M$.
We introduce a new start 
state, which, if the first letter seen lies in the alphabet of $M$,
then the automaton reacts is if it had been in the start state of 
the machine for $M$,
otherwise it reacts as if it had been in the start state of the 
machine for $L$.

 \qed

The interested reader can consult Hopcroft and Ullman's book \cite{HU}
for other properties of  (deterministic) pushdown automata and their 
languages.

\HH 4. $\ppp$ groups

Let $G$ be a group, and let $A$ be a finite monoid generating set, 
i.e.,
a finite set equipped with a map from $A\to G$, such that the induced
monoid homomorphism $A^*\to G$ is surjective.
The homomorphism  $A^*\to G$  is here denoted by $w \mapsto \eval w$.
A language $L \subset A^*$ is a \Em{parallel poly pushdown
structure} for $G$ if
\item {1.} the map $L\mapsto \eval L = G$ is a bijection;
\item{2.} $L$ is ${\cal P}$;
\item{3.} for each $a \in A$, $\{(w,w')\in L\times L \mid \eval {w'}
= \eval{wa} \}$ is AT$\ppp$. 

\noindent We say $G$ is \Em{parallel poly pushdown} (${\cal P}$) if $G$ has a
parallel poly pushdown structure.

\rk Remarks

1) Concerning  condition 3, notice that it is easy to check whether 
or not a given pair of words lies in $L\times L$.  
Thus, when checking this condition, we
will always assume we are given a pair in $L\times L$ and concentrate
on checking equality. 

2) If $G$ is a $\ppp$ group, and $B$ is any finite monoid generating set 
for $G$, we do not know whether or not there is a $\ppp$ structure for
$G$ with this generating set.
 
 3) Clearly the synchronous and asynchronous automatic groups of 
 \cite{ECHLPT} are $\ppp$ (with respect to all finite generating sets).
 
4) Without loss of generality, we can take $A$ to be closed under 
inverses.  We need only check that if $\{(w,w')\in L\times L \mid 
\eval {w'} = \eval{wa} \}$ is AT$\ppp$, then so is $\{(w,w')\in 
L\times L \mid \eval {w'} = \eval{wa^{-1}} \}$.  But this is exactly the same set of 
pairs taken in the opposite order, so it suffices to unplug the two 
input tapes and plug them back into each other's sockets!

\endrk

\th Proposition 4.1

If $G$ is $\ppp$ then $G$ has a $\ppp$ structure in which the
identity is represented by the empty word. \endth

\pf Proof

Suppose that $L'$ is a $\ppp$ structure for $G$ and that $w \in L'$ is
the representative for the identity.  We claim $L=L'-\{w\} \cup 
\{\epsilon\}$
is also a $\ppp$ structure for $G$.  Clearly $L$ is $\ppp$, by 
Proposition 3.1.  
Now suppose $L'_a=\{(w,w')\in L'\times L' \mid \eval {w'} = \eval{wa} \}$
and that for each $a \in A$, the representatives of $\eval a$ and
$\eval a^{-1}$ are $u_a$ and $v_a$.  Then using $L_a$ to represent
$\{(w,w')\in L\times L \mid \eval {w'} = \eval{wa} \}$, we have that
$L_a = L'_a -\{(w',u_a), (v_a,w')\} \cup \{(\epsilon,u_a),(v_a,\epsilon)\}$.
As before, it is easy to see that this is AT$\ppp$.\qed

\th  Theorem 4.2

Suppose $G$ is ${\cal P}$.  Then $G$ has a solvable word problem.
\endth 

\pf Proof

The proof of this is general nonsense.  Solvability of the word 
problem is independent of generating set.  We choose a generating 
set $A$ so that $L\subset A^{*}$ is a $\ppp$ structure.  We have now 
organized things so that

{\it 1.} $L$ is a recursively enumerable language which surjects to 
$G$.  (In fact, it is a recursive language which bijects to $G$.)

{\it 2.}  For each generator $a$ the subset of $L\times L$ which 
denote $a$-edges is recursive.

{\it 3.} $L_{1}$, the set of normal form words for the identity, is 
recursive.  (In fact, it is the language containing only 
the empty word.)

We are given the word $w=a_{1}\cdots a_{k}$ and asked to determine if 
$\eval w=1$.  We let $w_{0}$ be a normal form word for the identity.  
We assume that $j\le k$ and that we have found $w_{j-1}$, a normal 
form word for $\eval {a_{1}\cdots a_{j-1}}$.  If $j<k$, we enumerate 
the normal form words of $L$ and for each word $w\in L$ of this 
enumeration, we test whether or not $(w_{j-1},v)$ denotes an 
$a_{j}$-edge.  We will eventually find such a $v$, and when we do, we 
take it as $w_{j}$.  If $j=k$, we test whether $w_{j} \in L_{1}$.  As 
$\eval {w_{k}}=\eval w$ this determines whether or not $w$ evaluates 
to the identity.  \qed

The efficiency of this algorithm depends on the efficiency of the 
enumeration of {\it 1.}\ and the decision procedures of {\it 2.}\
and {\it 3.}\  For $\ppp$ structures, as with automatic structures, 
deciding {\it 1.}\ or {\it 2.}\ can be done in a length of time 
proportional to the length of the proposed word or words.  Likewise 
the decision entailed in {\it 3.} takes place in one step.  In the 
case of an automatic structure, producing $w_{j}$ from $w_{j-1}$ is 
highly efficient.  It can be done in linear time, and this gives rise 
to a quadratic time algorithm for solving the word problem in an 
automatic group.  This final estimate depends on the fact that the 
length of a normal form word is linear in the length of the element it 
represents.  

For an asynchronously automatic group, this algorithm can rise to 
exponential time for the simple reason that the length of $w_{j}$ may 
be exponential in $j$.  Since $\ppp$ groups include the asynchronously 
automatic groups, this algorithm can be at least that bad for $\ppp$ 
groups.  At present we do not know how to bound the length of a $\ppp$ 
normal form word for a group element of length $k$.  This would shed 
some light on the efficiency of this algorithm for $\ppp$ groups.

It follows from Theorem 4.2 that $\ppp$ groups are recursively 
presented.  We will see below ( Corollary 5.5) that they are not 
necessarily finitely presented.

\HH  5. Closure properties of $\ppp$ groups

\th  Theorem 5.1

The set of $\ppp$ groups is closed under direct product and free
product. \endth

\pf Proof

We suppose that $G$ and $G'$ are $\ppp$ groups and that $L \subset 
A^{*}$
and $L'\subset {A'}^{*}$ are $\ppp$ structures for $G$ and $G'$
respectively, with disjoint finite alphabets $A$, $A'$.  

We claim that $M=LL' \subset (A \cup A')^*$ is
a $\ppp$ structure for the direct product $G \times G'$.  
Map $A$ and $A'$ to the images of $\eval A$ and $\eval{A'}$ under
their natural inclusions into $G\times G'$. 
In Proposition 3.1 it is shown that  $M$ is
$\ppp$, and the induced map $M\to G \times G'$
is clearly bijective.  
It remains to show that, for each $a \in A \cup A'$, 
$\{(w,w')\in M\times M \mid \eval {w'} = \eval{wa} \}$ is \atp. 

Suppose $a \in A'$.
Since $L'$ is a $\ppp$ structure for $G'$ there is a collection of
asynchronous two tape pushdown automata which taken together define
$L''=\{(w,w')\in L'\times L' \mid \eval {w'} = \eval{wa} \}$.
Modify each of these machines to look initially for the diagonal in 
$A^* \times A^*$, followed by the representatives of the elements
of $L''$ when the first letter in $A'$ appears. 
A similar argument applies when $a \in A$.

Now consider the free product.
We can assume that $L$ and $L'$ contain the empty
word as their respective representatives for 1.  
The language 
$N=L  (L'^-L^-)^\star L'$ 
is $\ppp$, by Proposition 3.1, and it clearly maps bijectively to $G*G'$.  
As usual, we must check that for each
$a \in A \cup A'$, $\{(w,w')\in N\times N \mid \eval {w'} = \eval{wa}
\}$ is \atp.  
Suppose that $a \in A$.  
We can assume $\eval a \ne 1$ for
otherwise there is nothing to prove.  
Let $M_1,\ldots,
M_k$ be the asynchronous two tape pushdown machines that determine
the corresponding language for $G$.  
We describe machines $M_0, M'_1,\ldots, M_k$ that determine 
$\{(w,w')\in N\times N \mid \eval {w'} = \eval{wa} \}$

Now $w$ and $w'$ must have the form $w=u_1  v_1 \dots
 u_l  v_l$ and  $w'=u'_1   v'_1  \dots   u'_j 
v'_j$. 
Further, if $\eval {w'}=\eval{wa}$, we have $|l-j|\le 1$.  
We use an asynchronous two tape finite state automaton  
$M_0$ to check that, except for the last $G$ factor, 
we have $u_i=u'_i$ and $v_i=v'_i$, 
and that for the last $G$ factor $u_i\ne u'_i$.

Each $M'_i$ operates as follows.  It starts by reading the two
initial $L$ portions of $w$ and $w'$ emulating the action of $M_i$.
Upon encountering a change of generating set, it acts in the same way $M_i$
would react to $\$ $.  
If $M_i$ would accept the initial pair, it
remembers this and clears its stack.  
It then checks to see if it has
read the final $L$ factors, and accepts if it has done so.  
If it has not, it reads through the next $L'$ portion of $w$ and $w'$.  
Upon encountering a change of generating set, 
it repeats its emulation of $M_i$.  
Clearly $M_0$ and $M'_i$ both accept the pair $(w,w')$ 
if and only if these
two differ only in their last $L$ factors, these last factors would
have been accepted by $M_i$ and the remaining $L'$ factors are
trivial.  
Hence $M_0, M'_1,\ldots, M'_k$ perform as required.  
An analogous construction works for $ a' \in A'$.\qed

\th Theorem 5.2

Suppose $K$ is a finitely generated abelian group and 
$$1 \to K \to G \buildrel p \over\to H \to 1$$
is a central extension of an automatic group $H$.  Then $G$ is $\ppp$.
\endth 

\pf Proof

Recall that we can identify $G$ with the set $K \times H$ endowed with
the multiplication $(k,h)(k',h')=(k+k'+\rho(h,h'),hh')$ where $\rho$
is the $2$-cocycle that determines the extension.  
We write $K=\bbZ^n \times F$ where $F$ is finite.  
We take $A_K=\{x_1^{\pm 1},\ldots,x_n^{\pm 1}\} \cup A_F$ 
to be a generating set for $K$,
where $\{\eval{x_i}\}$ is a basis for $\bbZ^n$ and $A_F$ has the same 
number of elements as $F$.
We assume $A'_H$ is a generating set for $H$ and take $A_H$ to
be a set of lifts of these generators chosen by the inclusion 
$H\mapsto \{0\} \times H$.  
Since $H$ is assumed to be automatic, there
is a regular language $L'_H \subset {A'_H}^*$ 
which bijects to $H$ and has the \lq\lq fellow traveler property"
(see below).  
We let $L_H$ be the corresponding sublanguage of $A_H^*$.  
We take $L_K=\{x_1^{m_1}\ldots x_n^{m_n}f \mid
m_i \in {\Bbb Z}, f \in A_f\}$.  
We claim $L=L_HL_K\subset (A_H\cup A_K)^*$ is a $\ppp$ structure for $G$.

Clearly the natural map $L\to G$ is  a bijection.  
It is easy to see that for each $a \in
A_K$ $\{(w,w')\in L\times L \mid \eval {w'} = \eval{wa}\}$ is the 
intersection of $L\times L$ with the
language of a synchronous two tape finite state automaton, 
and thus \atp.  
We now show that for each $a \in A_H$, $\{(w,w')\in L\times L
\mid \eval {w'} = \eval{wa}\}$ is \atp.  
So suppose we have $w=ux_1^{m_1} \ldots x_n^{m_n} f$, 
$w'=u' x_1^{m'_1} \ldots x_n^{m'_n}f'$, 
and $a \in A_H$.  We will use $p_i$ and $p_f$ to represent
the maps  $K \to \langle \eval{x_i}\rangle$
and $K \to F$.  
We will have 
$$\eqalign{\eval {w'}&=\eval {wa}\hskip 1 true in\iff \cr
\eval{u' x_1^{m'_1} \ldots x_n^{m'_n}f'}&=
\eval{ux_1^{m_1} \ldots x_n^{m_n} fa}\hskip .09 true in\iff\cr
\eval{x_1^{m'_1} \ldots x_n^{m'_n}f'}&=
\eval{uau'^{-1}x_1^{m_1} \ldots x_n^{m_n} f .}}$$
These will happen if and only if:
\item{1.} $p(\eval{w'})= p(\eval{wa})$,
\item{2.} for each $i$, $m'_i=m_i + p_i(\eval{{ua}{{u'}^{-1}}})$, and
\item{3.} $f'=f+p_f(\eval{{ua}{{u'}^{-1}}})$.  

The first condition is easily checked by the comparator two tape
finite state automaton for $H$.  We assume the first condition is
satisfied.  To check the second condition, for each $i$ we build a two
tape pushdown automaton $M_i$ which reads $(u,u')$ and leaves
${p_i(\eval{{ua}{{u'}^{-1}}})}$ on its stack.
  
Recall that if condition 1 is satisfied, 
$L'_{H}$ has the fellow traveler property,
i.e., there is a constant $C$ so that for each
$j$, there is $t_j\in H$ with $\ell(t_j) \le C$,
so that $p(\eval{u'(j)}) = p(\eval{u(j)})t_j$.
(Here $u(j)$ and $u'(j)$ denote the initial segments of length $j$
 of $u$ and $u'$.)   
 
We start with the stack empty and $t_{0}=1$.  We assume inductively 
that $M_{i}$ has read the first $j$ letters of $u$ and $u'$, that it 
has the value ${p_i(p(\eval{u(j)})t_j p(\eval{u'(j)^{-1}})}$ on the 
stack and knows the value of $t_j$.  (Since $\ell(t_j)\le C$, this 
latter requires only a bounded amount of memory and can be carried in 
$M_{i}$'s internal memory.) $M_{i}$ reads $a_{j+1}$ and $a'_{j+1}$ 
from $u$ and $u'$ respectively.  If one of $u$ or $u'$ has been 
exhausted, but the other has not, then one of these letters can be 
taken to be the empty word.  Then $t_{j+1}=p(\eval{a_{j+1}})^{-1} t_j 
p(\eval{a'_{j+1}})$.  Since there are only finitely many values in 
this formula, this computation can be done in $M_{i}$'s internal 
memory.  Now observe that 
$$\eqalign {p_i(p(\eval{u(j+1)})t_{j+1}&p(\eval{u'(j+1)^{-1}}) =\cr 
&p_i(p(\eval{u(j)})t_jp(\eval{u'(j)^{-1}}) + p_i\rho(t_{j+1}, 
p(\eval{a_{j+1}'}^{-1})t_j^{-1}p(\eval{a_{j+1}})).\cr}$$
The first term on the right is the contents of the stack before 
$M_{i}$ reads $a_{j+1}$ and $a'_{j+1}$ and the second term on the 
right is determined by the finitely many possible values of $a_{j+1}$, 
$a'_{j+1}$, $t_{j}$ and $t_{j+1}$.  Thus, in order to compute 
$p_i(p(\eval{u(j+1)})t_{j+1}p(\eval{u'(j+1)^{-1}})$, $M_{i}$ need only 
add $p_i\rho(t_{j+1}, p(\eval{a_{j+1}'}^{-1}) t_j^{-1} 
p(\eval{a_{j+1}}))$.  This can be done since only finitely many such 
values occur.  If condition 1 is satisfied, then when $M_{i}$ is 
finished reading $u$ and $u'$, it's final $t$-value will be $p(a)$ and 
the stack will contain ${p_i(\eval{ua}\eval{{u'}^{-1}})}$ as required.

It is now easy to check condition 2 using $\{M_i\}$.  For each $i$ we
build a machine $M'_i$ which first uses $M_i$ to read $u$ and $u'$.
$M'_i$ reads $w$ and $w'$ until it gets to the $x_i^{\pm 1}$ portion.
It then pops the contents of the stack, canceling it letter by letter
against the $x_i^{\pm 1}$ letters of $w$ or $w'$ as appropriate.
After it has emptied the stack, it reads the $x_i^{\pm 1}$ letters of
$w$ and $w'$ one at a time and accepts if and only if these agree in
sign and it exhausts these simultaneously.  

Condition 3 can be checked in the same way, except that here we do
not require a stack to keep track of $p_f$ since only finitely
many values occur. \qed

\th Corollary 5.3

Every finitely generated nilpotent group of nilpotency class 2 is $\ppp$.
\endth

Notice that Theorem contrasts with the fact that 
the class of (asynchronously) automatic groups contains
no nilpotent groups of class greater than 1 (see \cite{ECHLPT}).

We now show how to construct some non--finitely presented $\ppp$ groups.

\th Theorem 5.4

The class of $\ppp$ groups is closed under wreath product.\endth

\pf Proof

Let $G$ and $H$ be $\ppp$ groups.  The wreath product $G \wr H$ can be 
identified with the semi--direct product $G^{[H]}\rtimes H$.  Let $L_H 
\subset A_H^*$ and $L_G\subset A_G^*$ be languages of $\ppp$ 
structures for these groups.  Fixing an ordering on the finite set 
$A_H$ induces a ``short lex" ordering $\prec$ on $A_H^*$, i.e., 
$u\prec v$ if $u$ is shorter than $v$ or if $u$ and $v$ have the same 
length, and $u$ is lexically prior to $v$.

To keep track of $H$--conjugates of elements of $G$, we introduce 
formal inverses of the elements of $L_{H}$.  That is, we introduce a 
disjoint copy ${A_{H}}^{\bf -1}$ of formal inverses of the elements of 
$A_{H}$, and the formal inverse map ${A_{H}}^{*} \to 
{({A_{H}}^{\bf -1})}^{*}: u=a_{1}\dots a_{n}\mapsto u^{\bf 
-1}={a_{n}}^{\bf -1}\dots {a_{1}}^{\bf -1}$.  We use the bold face 
exponent $^{\bf -1}$ to distinguish the process of taking formal 
inverses in ${A_{H}}^{\bf -1}$ from the process of taking inverses in 
$A_{H}$.  We will also need a second copy of $A_{H}$ which we will 
denote by $A_{H}^\circ$.  The isomorphism between  $A_{H}$ and 
$A_{H}^{\circ}$ induces an isomorphism of  $A_{H}^{*}$ and 
$(A_{H}^{\circ})^{*}$ denoted by $u \mapsto u^{\circ}$.

We now take

$$\eqalign{L=\{w_1 \ldots w_k u_0^{\circ} \mid & ~ w_i= u _i v_i{u_{i}}^{\bf -1}  \text{ 
for }i=1,\ldots k,\cr
& ~ u_i \in L_H \text{ for }i=0,\ldots k, \cr
& ~ v_i \in L_G \text{ for } i=1,\ldots ,k, \cr
& ~ u_i \prec u_{i+1}\text{ for }i=1,\ldots, k-1\}.\cr}$$

Because the natural maps from $L_H$ and $L_G$ to $H$ and $G$ are 
bijections, it is easy to see that the natural map from $L$ to $G \wr 
H$ is also a bijection.  Notice that since the $\{u_i\}$ are in 
strictly increasing $\prec$ order, $\eval{u_i} \ne 1$ for 
$i=2,\ldots,k$.

We now check that $L$ is a $\ppp$ language.  We use one machine to 
check that each $u_{i}$ and $u_{i}^{\bf -1}$ are formal inverses of 
each other.  This is done by successively pushing each $L_{H}$ subword 
onto the stack and popping it off letter by letter and checking it 
against the next $L_{H}^{\bf -1}$ subword encountered.  This machine 
accepts only if each of these comparisons is successful (except of 
course for the $u_{0}$ subword).

A second machine checks that for $i=1,\ldots , k-1$, we have $u_i 
\prec u_{i+1}$.  This machine acts by pushing each $u_{i}$ subword 
encountered onto the stack.  If the first machine accepts the word, 
then each $u_{i}^{\bf -1}$ subword is the formal inverse of the 
previous $u_{i}$ subword.  Since pushing this onto the stack reverses 
order, one may now pop the contents of the stack, comparing them 
letter by letter with the next $u_{i}^{\bf -1}$ subword to determine if 
they are in $\prec$ order.  This machine accepts if the each of these 
comparisons is successful.

Finally, a collection of $L_H$ machines successively check each 
$L_{H}$ (and $L_{H}^{\circ}$) subword and accepts only if each of 
these leads to an acceptance, while a collection of $L_G$ machines 
check each $L_G$ subword.  All of these machines accept if and only if 
$w \in L$.

We now wish to check that for each $a$, $\{(w,w')\in L\times L \mid 
\eval{w'}=\eval{wa}\}$ is \atp.  If $\eval a$ is in $\eval {A_H}^{\pm 
1}$, this is easy to do.  We simply check that $w$ and $w'$ are 
identical up to the final $u_{i}$ subword, and then use the 
multiplication machines for $L_H$ to compare $u_{o}^{\circ}$ with 
$u_{o}^{\prime\circ}$.  (When multiplying by $a^{\bf -1}\in A_{H}^{\bf 
-1}$, notice that $\eval{u'}=\eval{ua^{\bf -1}} \iff 
\eval{u}=\eval{u'a}$.) So we suppose that $\eval a \in \eval{A_G}$ and 
$w=w_1 \ldots w_k u_0^{\circ}\in L$, with 
$w_{i}=u_{i}v_{i}{u_{i}}^{\bf -1}$, and with the $u_{i}$ in increasing 
order.  We can assume that $\eval a \ne 1$.  Now if $\eval{w'}=
\eval{wa}$, there are three possibilities for $w$:

\item {1.}  there is no $i$, $1\le i \le k$ so that $u_i=u_0$;  

\hskip 0.1cm in this case,  
$w'=w_1 \ldots w_i u_0 v_a{u_0}^{\bf -1} w_{i+1} \ldots w_k u_0^{\circ}$, 
where $v_a$ is the $L_G$ word for 
$\eval a$ and $i$ is the largest value such that $u_i \prec u_0$;

\item{2.}  there is an $i$, $1 \le i \le k$ so that $u_0=u_i$ and
$\eval {v_i}= (\eval a)^{-1}$; 

\hskip 0.1cm in this case,  $w'= w_1 \ldots w_{i-1} w_{i+1} \ldots w_k 
u_0^{\circ}$;

\item{3.}  there is an $i$, $1 \le i \le k$ so that $u_0=u_i$ and
$\eval {v_i}\ne (\eval a)^{-1}$.  

\hskip 0.1cm in this case 
$w'= w_1 \ldots w_{i-1}u_iv'_i {u_i}^{\bf -1} w_{i+1} \ldots w_k 
u_0^{\circ}$,  
where $\eval{v'_i}= \eval {v_{i}a}$. 

\noindent In particular, 
\item{($\star$)} $w'$ differs from $w$ by either differing at $v_i$ with
$u_i=u_0$ or by insertion of $u_0 v_a {u_0}^{-1}$
or by deletion of $u_0 v_{a^{-1}} {u_0}^{-1}$
($\eval{v_a}=\eval a$ , $\eval{v_{a^{-1}}} = (\eval{a})^{-1}$).

We first build a machine $M_0$ that determines if $w$ and $w'$ differ 
in exactly one of these ways.  This machine starts reading each of $w$ 
and $w'$.  The disjoint sets of generators indicate whether a letter 
being read lies in a $u_{i}$, a $v_{i}$, a ${u_{i}^{\bf -1}}$ or 
$u_{0}^{\circ}$ subword.  The machine $M_0$ reads each of $w$ and $w'$ 
one letter at a time until it discovers a discrepancy.  If that first 
discrepancy is in the last $u_{i}$ subword, it enters a fail state.  
If the first discrepancy is in a $v_{i}$ subword, it finishes reading 
that subword from each side and accepts if the remainder of the two 
words is identical.  If $M_0$ discovers the first discrepancy in a 
$u_{i}$ subword (other than the last one), say, $u_i \ne u^{\prime 
}_i$, it continues to read the $u_{i}$ subwords and determines which 
word is $\prec$ earlier.  If $u_i \prec u^{\prime}_i$, then 
$M_{0}$ checks if $v_i$ was the word in $L_{G}$ which represents 
$(\eval a)^{-1}$.  If not, it rejects the pair.  This information is 
preserved in state of the machine while the subword $v_{i}$ was read.  
If $\eval{v_i}= (\eval a)^{-1}$, it then continues reading, pushing 
${{u_i}}^{\bf -1}$ onto the stack.  The subword $u_{i+1}$ is then read, 
and after this, $M_{0}$ continues reading each subword, checking for 
equality.  When it arrives at the final $u_{i}$ subword, it pops the 
stack checking to see if $u_i=u_0$ and accepts if both these things 
happen.  If, on the other hand $u^{\prime}_i\prec u_i$, $M_0$ checks 
that $v'_i$ is the word in $L_{G}$ representing the element $\eval a$.  
If so, it pushes ${u_i^{\prime }}^{-1}$ onto the stack and proceeds as 
before to check that $u'_i=u_0$ and that the remainder of the words 
are identical.  Finally, if $M_0$ does not encounter a discrepancy, it 
rejects the pair.  Thus if $M_0$ is given a pair $(w,w') \in L\times 
L$, it determines whether or not the pair satisfies ($\star$), above.

It is now easy to check if $\eval {w'}=\eval{wa}$.  
For each machine
$M_j$ used to check $a$ multiplication in $L_G$, 
we build a machine $M'_j$.  
$M'_j$ checks each $u_i$ against the corresponding
$u_i^{\prime}$ looking for inequality.  
It checks each $v_i$
against the corresponding $v'_i$ in the manner of $M_j$.  It continues
both of these tasks until it succeeds at one of them.  Notice that if
the pair is in $L\times L$ and is accepted by $M_0$, $M'_j$ can only
succeed at one of these tasks, and for only one value of $i$.  Thus
all of these machines together determine $\{(w,w')\in L\times L \mid
\eval{w'}=\eval{wa}\}$ as required. \qed

In view of the fact that there are wreath products of $\ppp$
groups which are not finitely presentable (for instance, if $C$ 
denotes the infinite cyclic group, $C\wr C$ is not finitely 
presented),
we have:

\th Corollary 5.5

There are $\ppp$ groups which are not finitely presented.\endth

We now turn to the consideration of semi--direct products of abelian 
groups by $\ppp$ groups.

\th Theorem 5.6

Suppose that $H$ is a $\ppp$ group, and $\varphi:H \to \Aut {\Bbb
Z}^n$.  Then ${\Bbb Z}^n \rtimes_\varphi H$ is $\ppp$. \endth

To obtain this result, we use the following lemma:

\th Lemma 5.7

Suppose $A\in \Aut {\Bbb Z}^n$.  Then 
$$L_A = \{(x_1^{m_1} \ldots x_n^{m_n}, x_1^{m'_1}\ldots x_n^{m'_n})
\mid (m'_1,\ldots,m'_n) = A (m_1,\ldots,m_n)\}$$
is \atp. 
\endth

\pf  Proof

We start by observing that for each $i$, $1\le i\le n$, we can build a 
pushdown automaton which reads $x_1^{m_1}\ldots x_n^{m_n}$ and ends 
with the $i^{\rm th}$ coordinate of $A (m_1,\ldots,m_n)$ on its stack.  
It does this by adding an $(A)_{ij}$ to the stack for each $x_j$ that 
it encounters.  Now, to recognize the $i^{\rm th}$ coordinate of 
the language $L_A$, we continue by reading 
$(x^{m'_1},\ldots,x^{m'_n})$ and subtracting $1$ for each $x_{i}$ or 
adding $1$ for each $x_{i}^{-1}$ we encounter.  Each machine accepts 
by empty stack.  The automaton then accepts by empty stack.  The 
collection of automata, one for each $i$, then determines the language 
$L_{A}$.\qed

\pf Proof of Theorem

Let $L_H\subset A_H^*$
be the language of a $\ppp$ structure for $H$, and take the language
$$L= L_H  \{x_1^{m_1}\ldots x_n^{m_n}\}.$$
Right multiplication by a generator  ${x_i}^{\pm 1}$ is easy
realized, and in view of the lemma, it is also easy to check right
multiplication by a generator $h \in A_H$ using $A=A_{\varphi(h)}$.
\qed

\th Corollary 5.8

The class of $\ppp$ groups contains groups of isoperimetric inequality
of every polynomial degree.\endth

\pf Proof

Bridson and Gersten \cite{BGe} have characterized the isomperimetric 
inequalities of groups ${\Bbb Z}^n \rtimes_A {\Bbb Z}$, where $A$ is a 
nilpotent matrix.  Such a group has isoperimetric inequality $n^{d+1}$ 
where $d$ is the size of the longest block of the Jordan canonical 
form of $A$.  \qed

In order to show that all geometric 3--manifold groups are in $\cal 
P$, we need the following result, which is much in the spirit of 
Theorem 5.6.

\th Lemma 5.9

If $G$ contains a finite index subgroup which is a semi--direct 
product of the form ${\Bbb Z}^n \rtimes_A {\Bbb Z}$,
 then $G$ is $\ppp$.
 \endth

\pf Proof

We take generators ${x_1}^{\pm 1}, \ldots ,{x_n}^{\pm 1}$ 
for ${\Bbb Z}^n$, $z^{\pm 1}$ for $\Bbb Z$ 
and a finite set $T$ giving a transversal for ${\Bbb Z}^n \rtimes_A
{\Bbb Z}$ in $G$.  
We take as our language $L=\{z^p x_1^{m_1}\dots x_n^{m_n}\}T$.
We must check that we can detect right multiplication by a generator.
So suppose we have the pair $w=z^r x_1^{m_1}\dots x_n^{m_n}t$ and
$w'=z^{r'} x_1^{m'_1}\dots x_n^{m'_n}t'$, and we wish to check if
$\eval {w'}=\eval{wg}$.  
To do this we will use $n$ automata, one for each of the letters 
$x_{i}$.  
Let $g$ be a generator.
There are only finitely many possibilities for $tg$, and 
each of these can be written in the form $us$ with
$u=\eval{z^ax_1^{b_1}\dots x_n^{b_n}} \in {\Bbb Z}^n \rtimes_A {\Bbb
Z}$ and $s \in T$.  
We start by reading the $z$ portions of $w$ and $w'$.  
If $|r -r'|$ exceeds the largest $z$ exponent in any $u$, then
$\eval{w'}\ne\eval{wg}$ and we reject the pair.  
If $|r'-r|$ does not exceed this value, we remember $r'-r$.  
This requires a bounded amount of memory (i.e., this information
is stored on one of a finite number of states).  
We now proceed to read the remainder of $w$, pushing
$x_i^{m_i}$ onto the stack in the $i^{\rm th}$ automaton.  
When we encounter $t$, we then compute $u$ and $s$.  
There are only finitely many possibilities for
these.  
We can now check if $r'-r=a$.   If not, we reject the pair.  
If $r'-r=a$, we continue.  
We read the ${\Bbb Z}^n$ portion of $w'$,
applying the transformation $A^a$ as in Lemma 5.6 popping and pushing  the
contents of each stack accordingly.  
When this is done, we accept if
and only if for each $i$, the stack of the $i$th machine 
contains $b_i$ and $t'=s$.
\qed

\th Corollary 5.10

Suppose $M$ is a 3-manifold which obeys the Thurston geometrization 
conjecture.  Then $\pi_1(M)$ is $\ppp$.
\endth

\pf Proof

If $M$ is such a manifold, $\pi_1(M)$ is the free product of an
automatic group with finitely many fundamental groups of closed Sol or
Nil geometry manifolds.  Each of these Sol or Nil geometry groups
contains a finite index subgroup of the form ${\Bbb Z}^2 \rtimes_A
{\Bbb Z}$, where $A$ is either nilpotent or Anosov.  Since automatic
groups are $\ppp$ and $\ppp$ groups are closed under free product, the
result now follows.\qed

\HH 6. Nilpotent groups

We here study the group  $U(n)$ whose elements are the $n\times n$
upper triangular integral matrices with $1$'s on the diagonal.  Our
interest in this group comes from the fact that if $G$ is a finitely
generated torsion free nilpotent group, then for sufficiently large
$n$, $G$ embeds in $U(n)$ \cite{Ba}.

\th Theorem 6.1

For each $n$, $U(n)$ is $\ppp$.  
In particular $\ppp$ contains
nilpotent groups of every class and every finitely generated torsion 
free nilpotent group embeds in a $\ppp$ group.  
\endth

\noindent Before proving the theorem, we recall some basic facts about $U(n)$.
Each element of $U(n)$ has the shape 
$$\pmatrix{1 &* &*  &\ldots &* &* \cr
           0 &1 &*  &\ldots &* &* \cr
           \vdots & &\ddots & &  &\vdots \cr
           0 &0 &0 &\ldots &1 &* \cr
           0 &0 &0 &\ldots &0 &1 \cr}$$
We will take $H(n)$ to be those matrices with the shape
$$\pmatrix{1 &* &* &\ldots &* &* &* \cr
           0 &1 &0 &\ldots &0 &0 &* \cr
           \vdots &  &\ddots &  & &  &\vdots \cr
           0 &0 &0 &\ldots &1 &0 &* \cr
           0 &0 &0 &\ldots &0 &1 &* \cr
           0 &0 &0 &\ldots &0 &0 &1 \cr } \ .$$
That is, $H(n)$ consists of those matrices of $U(n)$ for which all
nonzero entries are either on the diagonal, the first row, or the last
column.  
Let $T(n)\subset H(n)$   be  those matrices in which all nonzero
entries are restricted to the diagonal and the top row.  
Similarly, let $R(n)\subset H(n)$ be those matrices in which all nonzero entries are
restricted to the diagonal and the extreme right column.  
Let $e_{ij}$ be the matrix which has as single nonzero entry, 
a $1$ in the $ij$ position. 

\th Lemma 6.2

{\it 1.} $T(n)$ is a free abelian group of rank $n-1$ generated
by $x_j=1+e_{1j}$, $2\le j\le n$.

{\it 2.} $R(n)$ is a free abelian group of rank $n-1$ generated
by $y_i=1+e_{in}$, $1\le i\le n-1$.

{\it 3.} $H(n)=T(n)R(n)$.  Further $H(n)$ has the presentation
$$\eqalign{\langle x_2,\ldots, x_{n-1},y_2,\ldots y_{n-1},z \mid
[x_i,x_j]=1, &[y_i,y_j]=1, [x_i,y_j]=1, i\ne j \cr
&[x_i,y_i]=z, z\text{ is central} \rangle \ . \cr}$$
In particular,
$L_{H(n)}=\{x_2^{p_2} \ldots x_{n-1}^{p_{n-1}} y_2^{q_2}\ldots
y_{n-1}^{q_{n-2}} z^r\}$
is the language of a $\ppp$ structure for $H(n)$.

{\it 4.} $U(n)$ is generated by $\{1+e_{ij} \mid j>i\}$. There is
a split short exact sequence 
$$1 \to H(n)\to U(n)\to U(n-2)\to 1 .$$
The splitting is given by the inclusion of $U(n-2)$ into $U(n)$ as
$1\oplus U(n-2)\oplus 1$, the subset of $U(n)$ for which the
off-diagonal entries in the first row and last column are all zero.
The action of the generator $1+e_{ij}$, $1<i<j<n$ on $H(n)$ carries
$x_2^{p_2} \ldots x_{n-1}^{p_{n-1}} y_2^{q_2}\ldots y_{n-1}^{q_{n-1}}
z^r$ to $x_2^{p'_2} \ldots x_{n-1}^{p'_{n-1}} y_2^{q'_2}\ldots
y_{n-1}^{q'_{n-1}} z^r$ where for $m \ne j$, $p'_m=p_m$,
$p'_j=p_i+p_j$, for $m\ne i$, $q'_m=q'$, $q'_i=q_i-q_j$.\endth

\pf Proof

First observe that 
$$e_{ij}e_{kl}=\cases{0 &if $j\ne k$,\cr e_{il} &if $j=k$,\cr}$$
so that for $i<j$, $(1+e_{ij})^{-1}=1-e_{ij}$.  Statements {\it 1.}\
and {\it 2.}\ now follow easily.  In this way it is also easy to see
that $x_2,\dots,x_{n-1},x_n=z=y_1,y_2,\dots,y_{n-1}$ fulfill the
relations of the presentation of {\it 3.}  One checks that $H(n)$ is
in fact a subgroup of $U(n)$.  From the relations of the
presentation it is easy to see that each product of the elements
$x_2,\dots,x_{n-1},x_n=z=y_1,y_2,\dots,y_{n-1}$ can be put into the
form 
$$x_2^{p_2} \ldots x_{n-1}^{p_{n-1}} y_2^{q_2}\ldots
y_{n-1}^{q_{n-2}} z^r.$$
A computation shows that this element corresponds to the matrix 
$$\pmatrix{1      &p_2 &p_3     &\dots &p_{n-1} &r+\sum_{i=2}^{n-1} p_iq_i\cr
           0      &1   &0       &\dots &0       &q_2 \cr
           \vdots &    &\ddots  &      &        &\vdots\cr
           0      &0   &0       &\dots &1       &q_{n-1}\cr
           0      &0   &0       &\dots &0       &1\cr}.$$
Thus $L_{H(n)}$ bijects to $H(n)$.  In particular, it follows that
the elements $x_2,\dots,x_{n-1},x_n=z=y_1,y_2,\dots,y_{n-1}$ generate
$H(n)$. Notice that the words of $L_{H(n)}$ are already in the form
$T(n)R(n)$, so $H(n)=T(n)R(n)$ as claimed.

Let ${\cal H}(n)$ be the group determined by the presentation.  Since
the generators of $H(n)$ obey the relations of the presentation, 
the obvious map from ${\cal H}(n)$ to $H(n)$ is a surjective 
homomorphism.
To see that this is an isomorphism, let $g$ be an
element of the kernel.  
Using the relations of the presentation, we
can write $g$ in the normal form we have used in $H(n)$.  
The fact that $g$ maps to the identity matrix in $H(n)$ forces
$p_2=\dots=p_{n-1}=q_2=\dots=q_{n-1}=0$ and $r=0$, 
so in fact $g=1$ in ${\cal H}(n)$ as required.

We must check that $L_{H(n)}$ gives a $\ppp$ structure.  
It is a regular language such that the natural map to $H(n)$
is bijective.  
We now show how to check right multiplication in the appropriate 
way.  
Right multiplying by $z$ or $y_i$ only increases $r$
or $q_i$ by $1$, and this can be checked by a finite state automaton.
It remains to check right multiplication by $x_i$.  
This increases $p_i$ by one and changes $r$ by $q_i$.  
This can be checked using a single stack.  
We leave the details to the reader.

\def\ts{\thinspace}

To  prove {\it 4.}, first notice that
$$\displaylines{
  \pmatrix{1      &a_{12} &a_{13} &\dots &a_{1\ts n-2} &a_{1\ts n-1}   &a_{1\ts n}\cr
           0 &1 &a_{23} &\dots &a_{2\ts n-2} &a_{2\ts n-1} &a_{2\ts n}\cr
           \vdots &       &\ddots &      &          &               &\vdots\cr
           0      &0      &0      &\dots &1         &a_{n-2\ts n-1} &a_{n-2\ts n}\cr
           0      &0      &0      &\dots &0         &1              &a_{n-1\ts n}\cr
           0      &0      &0      &\dots &0         &0           &1\cr} 
=\cr
  \pmatrix{1      &0      &0      &\dots &0         &0           &0\cr
           0      &1      &a_{23} &\dots &a_{2\ts n-2} &a_{2\ts n-1}   &0\cr
           \vdots &       &\ddots &      &          &            &\vdots\cr
           0      &0      &0      &\dots &1         &a_{n-2\ts n-1} &0\cr
           0      &0      &0      &\dots &0         &1           &0\cr
           0      &0      &0      &\dots &0         &0           &1\cr} 
  \pmatrix{1      &a_{12} &a_{13} &\dots &a_{1\ts n-2} &a_{1\ts n-1}   &a_{1\ts n}\cr
           0      &1      &0      &\dots &0         &0           &a_{2n}\cr
           \vdots &       &\ddots &      &          &            &\vdots\cr
           0      &0      &0      &\dots &1         &0           &a_{n-2\ts n}\cr
           0      &0      &0      &\dots &0         &1           &a_{n-1\ts n}\cr
           0      &0      &0      &\dots &0         &0           &1\cr}
.\cr} $$
We now include $U(n-2)$ into $U(n)$ as $1\oplus U(n-2)\oplus 1$.
Then the above equation can be read as saying that $U(n)=U(n-2)H(n)$.
Now $U(1)$, $U(2)$, and $H(n)$ (for all $n$) are generated by elements
of the form $1+e_{ij}$, $i<j$.  
Inductively, we conclude that $U(n)$ is generated by $\{1+e_{ij}\mid i<j\}$.  
It is now easy to check that $H(n)$ is a normal subgroup, since we need
only check conjugation by generators, $1+e_{ij}$, $i<j$.  
This will verify that the action on $H(n)$ is as claimed.  
It remains only to check that the quotient of $U(n)$ by $H(n)$ is as claimed.  
We take the map of $U(n)$ to $U(n-2)$ to be the ``forgetful'' map which
strips each matrix of its first and last rows and columns.  
A computation shows that this is a homomorphism and the kernel is
clearly $H(n)$.\qed

\pf Proof of the Theorem

We use the short exact sequence to perform an induction argument.  
Since this induction takes us from $U(n-2)$ to $U(n)$, it requires two
basis steps, $U(1)$ and $U(2)$.  
These groups are respectively the trivial group and $\Bbb Z$,
both of which are $\ppp$, so the basis step is complete. 

We now assume that $L(n-2)$ is a regular language over an alphabet
consisting of letters for $\{1\pm e_{ij}\mid 1 < i<j <n-2\}$ which is $\ppp$
structure for $U(n-2)$.  
We take $L(n)=L(n-2)L_{H(n)}$.  
By the short exact sequence and our choice of generators for $H(n)$, 
this is a language over the desired generating set   
and is in one--to--one correspondence with $U(n)$.  
We must check that this is a $\ppp$ structure.  
Right multiplication by an element of $H(n)$ is easily checked, 
as we have seen in {\it 3.}\ of the Lemma 6.2.  
On the other hand, we can check right multiplication by an element of $U(n-2)$ 
since this requires that we check right multiplication in $U(n-2)$ 
(which we can do by induction)
and check the action of a generator on $L_{H(n)}$, which we can do by
Lemma 5.7.\qed

\Bib References

\rf {BGSS}  G. Baumslag, S.~M. Gersten, M. Shapiro and H. Short,
{\it Automatic groups and amalgams}, J. of Pure and Applied Algebra 76, 
(1991), pp 229--316.

\rf {Ba} G.~Baumslag, {\it  Lecture Notes on nilpotent groups\/}, 
in Regional  Conf. series of  A.M.S. 2 (1971).

\rf {Br} N. Brady, Thesis, University of California at Berkeley, 1993.\NR

\rf {BGe} M.~R. Bridson and S.~M. Gersten, {\it The optimal isoperimetric
inequality for torus bundles over the circle\/}, to appear in Quart. J. Math..

\rf {BGi} M.~R. Bridson and R. Gilman, {\it Formal language theory and the 
geometry of 3--manifolds\/}, to appear in Comm. Math. Helv.

\rf {ECHLPT} D.~B.~A. Epstein, with J.~W. Cannon, D.~F. Holt, 
S.~V.~F.~Levy, M.~S. Paterson  and W.~P. Thurston,
``Word processing in Groups'',
 Jones and Barlett, Boston--London, 1992.

\rf {GB2} G.~Baumslag, {\em Topics in Combinatorial Group Theory},
Birkh\"auser, 1993.

\rf {HU} J.~E. Hopcroft and J.~D. Ullman, 
{\em Introduction to  automata theory, languages and computation},
Addison--Wesley, 1979.

\rf {S} M. Shapiro, unpublished.

\endBib

\Coordinates

City College of New York
New York, NY 
USA
\bigskip

University of Melbourne
Parkville, VIC
Australia
\bigskip

LATP, Centre de Math\'ematiques et d'Informatique,
Rue Joliot--Curie, Universit\'e de Provence, F--13453
Marseille cedex 13, France

\endCoordinates
\bye